\documentclass[11pt,a4paper]{article}

\usepackage{amsmath}
\usepackage{amssymb}
\usepackage{amsthm}
\usepackage{fancyhdr}
\usepackage[hidelinks,bookmarksnumbered]{hyperref}
\hypersetup{pdftitle={Needlepoints in Quasi-Normed Orlicz Spaces and a Question of Kalton and Peck},
  pdfauthor={Dong Ji and Yunan Cui}}
\numberwithin{equation}{section}
\makeatletter
\renewcommand\section{\@startsection{section}{1}{\z@}%
  {-3.25ex\@plus -1ex \@minus -.2ex}%
  {2.2ex\@plus .2ex}%
  {\centering\normalfont\scshape}}
\renewcommand\subsection{\@startsection{subsection}{2}{\z@}%
  {-2.75ex\@plus -1ex \@minus -.2ex}%
  {1.5ex \@plus .2ex}%
  {\normalfont\scshape}}
\makeatother

\theoremstyle{plain}
\newtheorem{theorem}{Theorem}[section]
\newtheorem{lemma}[theorem]{Lemma}
\newtheorem{proposition}[theorem]{Proposition}

\theoremstyle{definition}
\newtheorem{definition}[theorem]{Definition}

\newcommand{\blfootnote}[1]{%
  \begingroup
  \renewcommand{\thefootnote}{}%
  \footnote{\hspace*{-\footnotesep}#1}%
  \addtocounter{footnote}{-1}%
  \endgroup}

\newcommand{\E}{\mathbb{E}}
\newcommand{\Prob}{\mathbb{P}}
\newcommand{\co}{\mathrm{co}}
\newcommand{\closedco}{\overline{\mathrm{co}}}

\begin{document}
\thispagestyle{plain}

\begin{center}
  {\scshape\large Needlepoints in Quasi-Normed Orlicz Spaces and a\\[0.3ex]
  Question of Kalton and Peck}\\[2ex]
  {\scshape Dong Ji and Yunan Cui}
\end{center}

\begin{quote}
\small
\noindent {\scshape Abstract.} We prove that the approachable points and the
needlepoints coincide in every quasi-normed Orlicz function space on the
probability interval. Their common set is either $\{0\}$ or $E_\Phi$. The proof gives a characterization by finite probability laws and
also applies to Orlicz function spaces with the modular topology, without
assuming $\Delta_2$. We then construct a separable real quasi-Banach
space whose approachable points form a line and whose only needlepoint is zero.
Every compact convex subset of this space is relatively locally convex. The
construction gives a negative answer to the question of Kalton and Peck for
general quasi-Banach spaces.
\end{quote}

\section{Introduction}\label{sec:intro}

The Krein--Milman theorem asserts that, in a locally convex Hausdorff space,
every nonempty compact convex set is the closed convex hull of its extreme
points \cite{KM40}. It is therefore natural to ask what becomes of compact convex
sets once local convexity is lost. The answer is governed by the geometry of
small vectors.

In some quasi-Banach spaces, a nonzero vector can be approximated by the
arithmetic means of finite families of vectors whose quasi-norms are
arbitrarily small. The notions of approachable point and needlepoint
distinguish two ways in which this can happen. For approachability one
controls the absolutely convex hull of each finite family by a bound
independent of the accuracy of approximation. For a needlepoint one can
choose the entire convex hull arbitrarily close to the segment from zero to
the limiting vector. The second requirement implies the first. Kalton and
Peck asked whether the converse holds for a point of a quasi-Banach space
\cite[Definitions 2.3 and 3.1, Section 5(2), p.\ 101]{KP80}. We prove the
converse in all quasi-normed Orlicz function spaces on the probability
interval and construct a quasi-Banach space in which it fails.

Roberts constructed an $F$-space containing a compact convex set without
extreme points \cite{Ro77}. In the classical example $L^p(0,1)$, $0<p<1$, every
point is a needlepoint \cite[Proposition 5.3]{Ka03}. Theorem~\ref{thm:lux} identifies the
subspace that replaces $L^p(0,1)$ for a general $b$-Orlicz function.

\blfootnote{\textit{2020 Mathematics Subject Classification.} Primary
46E30; Secondary 46A16, 46A55, 60G50.\\
\textit{Key words and phrases.} $b$-Orlicz function, Luxemburg quasi-norm,
approachable point, needlepoint, compact convex set.\\
\textit{Corresponding author:} Yunan Cui.}

All vector spaces are real and functions are identified up to equality
almost everywhere. We use Lebesgue measure on $(0,1)$, and write
$\mathcal A(X)$ and $\mathcal N(X)$ for the approachable points and
needlepoints of the whole space $X$. The Orlicz functions $\Phi$ and
$\varphi$, the modular $I_\Phi$, the spaces $L^\Phi$ and $E_\Phi$, and the
Luxemburg quasi-norm $\|\cdot\|_\Phi$ are defined in Section~\ref{sec:prelim}; the
$b$-condition in Definition~\ref{def:borlicz} guarantees that the Luxemburg functional
is a complete quasi-norm \cite[Theorem~1]{JC24}.

\begin{theorem}\label{thm:lux}
Let $\Phi$ satisfy Definition~\ref{def:borlicz} with the dilation inequality required
only at infinity, that is, for some $b>1$ and $u_0>0$,
\begin{equation}\label{eq:tail-b}
  \Phi(bu)\ge2\Phi(u)\qquad(u\ge u_0).
\end{equation}
Let $X=(L^\Phi(0,1),\|\cdot\|_\Phi)$. Then
\[
  \mathcal A(X)=\mathcal N(X)=\mathcal N(E_\Phi)=\mathcal A(E_\Phi)
  \in\bigl\{\{0\},E_\Phi\bigr\},
\]
where $E_\Phi$ has the relative topology. If the common set is $\{0\}$, every
compact convex subset of $X$ is affinely homeomorphic to a compact convex set
in a locally convex space. In the other case $E_\Phi$ contains a nonempty
compact convex set without extreme points.
\end{theorem}

On the probability interval, $E_\Phi=L^\Phi$ precisely when
$\Phi\in\Delta_2$ at infinity \cite[Theorem~2.14]{CHKK19}; Proposition~\ref{prop:proper-example} gives a $b$-Orlicz function
for which $E_\Phi$ is proper. In the zero case the compact-set conclusion
follows from Kalton's embedding theorem \cite{Ka80}, and in the other case the
witnesses lie in the complete subspace $E_\Phi$, where the Roberts
theorem applies \cite{KPR84}.

The proof gives a statement about the modular topology itself. Here we use
the wider convention of Kalton \cite[pp.\ 126--127]{Ka77} for a monotone Orlicz
function $\varphi$ (Definition~\ref{def:monotone}), define $L^\varphi$ and $E_\varphi$ by
the same scale conditions as in Section~\ref{sec:prelim}, and give $L^\varphi$ the topology in
which $x_j\to 0$ when $I_\varphi(cx_j)\to 0$ for every $c>0$.

\begin{theorem}\label{thm:modular}
For every such $\varphi$, the conclusions of Theorem~\ref{thm:lux} hold with
$X=L^\varphi(0,1)$ in its modular topology and with $E_\varphi$ in place of
$E_\Phi$.
\end{theorem}

The proof has two function-space ingredients. A nonzero approachable
point yields finite probability laws together with a uniform moment bound,
and this bound places all independent copies in a single bounded convex
hull, giving the segment approximation required for a needlepoint. The
argument works in the modular topology itself, even when
$\varphi\notin\Delta_2$.

\begin{theorem}\label{thm:counter}
There are a separable real $1/2$-Banach space $X$ and a vector $u\in X$ of
quasi-norm one such that
\[
  \mathcal A(X)=\mathbb R u,\qquad \mathcal N(X)=\{0\}.
\]
Every compact convex subset of $X$ is locally convex in its relative
topology. In particular, every nonempty compact convex subset of $X$ is the
closed convex hull of its extreme points.
\end{theorem}

Here relative local convexity means that each point has a base of convex
neighborhoods in the relative topology. In the construction, the standard coordinate vectors
from a sequence of finite-dimensional spaces become arbitrarily small while
their means share one image in a quotient; a uniform estimate excludes
nonzero needlepoints and, on compact sets, gives a base of convex relative
neighborhoods. Thus every compact convex subset is relatively locally
convex, although $\mathcal A(X)$ is a nonzero line.

\section{Preliminaries}\label{sec:prelim}

Let $\mu$ be Lebesgue probability measure on $(0,1)$.

\begin{definition}\label{def:monotone}
A function $\varphi:[0,\infty)\to[0,\infty)$ is called a \emph{monotone
Orlicz function} if it is finite, nondecreasing, continuous at zero,
vanishes at zero, and is not identically zero. Throughout this section
$\varphi$ satisfies the hypotheses of Theorem~\ref{thm:modular}.
\end{definition}

\begin{definition}\label{def:borlicz}
A function $\Phi:\mathbb R\to[0,\infty)$ is called a \emph{$b$-Orlicz function}
if it is finite, continuous and even, vanishes at zero, is positive and
nondecreasing on $(0,\infty)$, and satisfies, for some $b>1$,
\[
  \Phi(u/b)\leq \tfrac12\Phi(u)\qquad (u\geq 0).
\]
\end{definition}

For any monotone Orlicz function $\Phi$ write
\[
  I_\Phi(x)=\int_0^1\Phi(|x(t)|)\,dt
\]
and set
\begin{align*}
  L^\Phi&=\{x:I_\Phi(cx)<\infty \text{ for some } c>0\},\\
  E_\Phi&=\{x:I_\Phi(cx)<\infty \text{ for every } c>0\}.
\end{align*}
\begin{definition}\label{def:indices}
The \emph{Matuszewska--Orlicz indices at infinity} of $\Phi$ are
\begin{align*}
  \alpha_\Phi&=\sup\bigl\{r>0:\exists\,C>0,\ \exists\,u_0\ \text{such that }\\[-2pt]
  &\hspace{2em}\Phi(au)\ge C a^r\Phi(u)\ \text{for all }a\ge1,\ u\ge u_0\bigr\},\\
  \beta_\Phi&=\inf\bigl\{r>0:\exists\,C>0,\ \exists\,u_0\ \text{such that }\\[-2pt]
  &\hspace{2em}\Phi(au)\le C a^r\Phi(u)\ \text{for all }a\ge1,\ u\ge u_0\bigr\},
\end{align*}
with the infimum equal to $\infty$ when the set is empty.
\end{definition}

For a $b$-Orlicz function $\Phi$, the Luxemburg functional on $L^\Phi$
\[
  \|x\|_\Phi=\inf\{a>0:I_\Phi(x/a)\leq 1\}
\]
is a complete quasi-norm \cite[Theorem~1]{JC24}. On the probability interval,
\cite[Theorem~1.8]{KZ19} gives the criterion $\alpha_\Phi>0$ for continuous
strictly increasing generators tending to infinity. A positive lower index
at infinity is equivalent to~\eqref{eq:tail-b}: iterating that condition
gives $\Phi(au)\ge\tfrac12 a^{\ln 2/\ln b}\Phi(u)$ for $a\ge1$ and
sufficiently large $u$, and a lower exponent $r>0$ gives
\eqref{eq:tail-b} upon choosing $b$ with $Cb^r\ge2$.
Definition~\ref{def:borlicz} requires the inequality for all $u\ge0$.

\begin{lemma}\label{lem:representative}
There is a continuous monotone Orlicz function $\Phi$ which equals $u$ on
$[0,1]$, is positive on $(0,\infty)$ with $\Phi(1)=1$, and defines the
same two spaces and the same topology as $\varphi$.
\end{lemma}

\emph{Proof.} If the given function already has these properties, we take
it directly. In general, use the averaging formula of Kalton
\cite[Lemma~3.2]{Ka77}, setting $G(0)=0$ and
\[
  G(u)=\int_1^2\varphi(su)\,ds=u^{-1}\int_u^{2u}\varphi(v)\,dv
  \qquad (u>0).
\]
The latter expression gives continuity away from zero. At zero use
$\varphi(u)\leq G(u)\leq\varphi(2u)$. These inequalities also show that $G$
defines the same two spaces and the same convergence as $\varphi$. Choose
$b>0$ with $G(b)>0$, put $H(u)=G(bu)/G(b)$, and let
\[
  \Phi(u)=
  \begin{cases}
    u, & 0\leq u\leq 1,\\
    H(u), & u>1.
  \end{cases}
\]
Fixed positive changes in the argument and value preserve the two spaces
and the topology. The change from $H$ to $\Phi$ preserves finiteness of each
scale integral, since their difference is bounded and $\mu((0,1))=1$. For
either $J=H$ or $J=\Phi$, convergence at every scale implies convergence in
measure, by $\mu\{|x_j|\geq r\}\leq I_J(x_j/r)$ for $r>0$. For fixed $a>0$
the integrands $H(a|x_j|)$ and $\Phi(a|x_j|)$ agree on $\{a|x_j|>1\}$. On its
complement their difference is bounded and tends to zero in measure. Its
integral therefore tends to zero. The two topologies agree. In particular
$\Phi$ is continuous, positive on $(0,\infty)$ and satisfies $\Phi(1)=1$.
\qed

\begin{definition}\label{def:fnorm}
The Mazur--Orlicz $F$-norm \cite{MO58}, defined for a monotone Orlicz function
in the integral form of \cite[p.\ 466]{No92}, is
\[
  \|x\|_F=\inf\{\lambda>0:I_\Phi(x/\lambda)\leq\lambda\}.
\]
\end{definition}

\begin{proposition}[see \cite{CHKK19}]\label{prop:fnorm}
The functional $\|\cdot\|_F$ is an $F$-norm on $X$: it is definite, monotone in
the order of the lattice, and satisfies the triangle inequality. Scalar
multiplication is jointly continuous, and
\begin{equation}
  \|x_j\|_F\to 0\iff I_\Phi(a x_j)\to 0\quad\text{for every }a>0.
  \label{eq:all-scale}
\end{equation}
\end{proposition}

\begin{proposition}\label{prop:complete}
The Orlicz space $(L^\Phi,\|\cdot\|_F)$ is an $F$-space. A set $B\subset L^\Phi$ is topologically
bounded if and only if $\sup_{v\in B}\|tv\|_F\to 0$ as $t\to 0$.
\end{proposition}

\emph{Proof.} A $\|\cdot\|_F$-Cauchy sequence $(x_j)$ is Cauchy in measure. Take
an almost everywhere convergent subsequence with limit $x$. For $r>0$ and
all sufficiently large $j,k$ we have $I_\Phi((x_j-x_k)/r)\leq r$. Fatou's
lemma along the subsequence gives $I_\Phi((x_j-x)/r)\leq r$ for the same
$j$. Thus $x\in L^\Phi$ and $x_j\to x$, proving completeness. The stated
characterization of topological boundedness is standard for $F$-norms. \qed

\begin{lemma}\label{lem:gauge}
Let $\varphi:\mathbb R\to[0,\infty)$ be even, continuous, positive on $(0,\infty)$
and nondecreasing there, with $\varphi(0)=0$. Its Luxemburg functional
\[
 \|x\|_\varphi=\inf\{a>0:I_\varphi(x/a)\le1\}
\]
is a quasi-norm if and only if $\varphi$ satisfies~\eqref{eq:tail-b}.
In that case it is complete and defines the same topology as $\|\cdot\|_F$.
\end{lemma}

\emph{Proof.} Assume~\eqref{eq:tail-b} and take $R>0$ so large that $\varphi(bu)\ge2\varphi(u)$ for
$u\ge R$, and define an even $\Psi$ by setting
$\Psi(t)=t\varphi(R)/R$ for $0\le t\le R$ and $\Psi(t)=\varphi(t)$
for $t\ge R$. It satisfies Definition~\ref{def:borlicz} with $B=2b$. For $t\ge R$
this follows from the tail inequality, and for $Bt\le R$ it follows
from linearity. In the remaining case $t<R<Bt$, use
$\Psi(Bt)\ge\varphi(R)$ when $t\le R/2$ and
$\Psi(Bt)\ge\varphi(bR)\ge2\varphi(R)$ when $t>R/2$.
Thus $\|\cdot\|_\Psi$ is a complete quasi-norm by
\cite[Theorem~1]{JC24}.

Choose $M\ge b$ so large that $\varphi(bR/M),\Psi(bR/M)\le1/2$.
For either ordering $(J,K)$ of $(\varphi,\Psi)$ and for $t<bR$ we have
$J(t/M)\le1/2$. For $t\ge bR$,
\[
 J(t/M)\le J(t/b)=\varphi(t/b)\le\tfrac12\varphi(t)=\tfrac12 K(t).
\]
Since the measure is one, $I_K(x/a)\le1$ implies $I_J(x/(Ma))\le1$.
Taking infima in the two directions gives
$M^{-1}\|x\|_\Psi\le\|x\|_\varphi\le M\|x\|_\Psi$.
The functionals have the same domain because the two integrands differ by
a bounded amount on $[0,R]$ and agree above $R$. Hence the original
functional is a complete quasi-norm.

For $\|x_j\|_\Psi\to0$, choose $a_j>\|x_j\|_\Psi$ tending to zero.
For fixed $c>0$ and $m\ge1$, eventually $c a_j\le B^{-m}$, so
$I_\Psi(cx_j)\le2^{-m}I_\Psi(x_j/a_j)\le2^{-m}$. Conversely,
convergence at every scale gives $\|x_j\|_\Psi\le r$ eventually for each
$r>0$ by testing at $1/r$. The bounded change below $R$ preserves
modular convergence at every scale by the argument of Lemma~\ref{lem:representative}:
either modular topology implies convergence in measure, as both functions
are positive at $R$. Hence the topology agrees with that of
$\varphi$, and thus with $\|\cdot\|_F$ by~\eqref{eq:all-scale}.

Conversely, suppose the functional is a quasi-norm with quasi-triangle
constant $K\ge1$. The function $\varphi$ must be unbounded, since
otherwise the functional vanishes on indicators of sets of sufficiently
small positive measure. Put $v(s)=\max\{t\ge0:\varphi(t)\le s\}$ for
$s>0$. For every set $A$ of positive measure,
$\|\chi_A\|_\varphi=1/v(1/\mu(A))$. Apply the quasi-triangle inequality
to the indicators of two disjoint sets of measure $1/(2s)$, $s\ge1$.
It gives $v(2s)\le 2K v(s)$. If $\varphi(t)>1$, take
$1\le s<\varphi(t)$. Then $v(s)\le t$, so $v(2s)\le2Kt$.
For $b>2K$ we obtain $\varphi(bt)>2s$. Letting $s\uparrow\varphi(t)$
gives~\eqref{eq:tail-b} for all sufficiently large $t$. \qed

For a real metrizable topological vector space $Y$ (in particular, an
$F$-space or a quasi-normed space), write $|y|$ for the distance of $y$
from zero in a defining translation-invariant metric. We recall the definitions of
\cite[Section 2]{Ka80}, \cite[Definitions 2.3 and 3.1]{KP80} and \cite[pp.\ 200--201]{KPR84}.

\begin{definition}\label{def:approach}
A point $y\in Y$ is \emph{approachable} if there is a
topologically bounded set $B\subset Y$ such that, for every $\varepsilon>0$,
there exist $y_1,\ldots,y_n$ with
\[
  \max_i|y_i|<\varepsilon,\qquad
  \left|y-\frac1n\sum_i y_i\right|<\varepsilon,\qquad
  \left\{\sum_i a_i y_i:\sum_i|a_i|\leq 1\right\}\subset B.
\]
\end{definition}

\begin{definition}\label{def:needle}
A point $y$ is a \emph{needlepoint} if, for every $\varepsilon>0$, there is
such a finite family with its entries and its mean error smaller than
$\varepsilon$ and with
\[
  \co\{y_1,\ldots,y_n\}\subset[0,y]+\{z:|z|<\varepsilon\}.
\]
Here $[0,y]=\{ty:0\leq t\leq 1\}$.
\end{definition}

We denote the set of approachable points and the set of needlepoints by
$\mathcal A(Y)$ and $\mathcal N(Y)$, respectively.

\begin{proposition}\label{prop:maps}
Continuous linear maps preserve approachability and the needlepoint
property. The set $\mathcal N(Y)$ is closed, and
$\mathcal N(Y)\subset\mathcal A(Y)$.
\end{proposition}

\emph{Proof.} One applies the map to the rows and, for approachability,
also to the common bounded set. For closedness, if $y_k\to y$, then
$t(y_k-y)\to 0$ uniformly in $0\leq t\leq 1$ by continuity of scalar
multiplication and compactness of $[0,1]$. A sufficiently accurate family
for $y_k$ therefore gives the two approximation conditions and the segment
condition for $y$. Moreover $\mathcal N(Y)\subset\mathcal A(Y)$. To see
this, choose needlepoint families for tolerances tending to zero, and in
each family form the convex hull $Q_m$ of its members together with zero.
Every point of $Q_m$ tends uniformly in distance to the segment $[0,y]$, in
the sense that $Q_m\subset[0,y]+U$ eventually for every neighborhood $U$ of
zero. This follows by the same uniform continuity for scalars in $[0,1]$.
Each $Q_m$ is compact, and $\bigcup_m Q_m$ is totally bounded since the
segment is compact. Hence $Q=\bigcup_m Q_m$ is topologically bounded,
as is $Q-Q$. Every
absolutely convex hull of the chosen families lies in $Q-Q$, as is seen by
separating the positive and negative coefficients. \qed

\section{A probabilistic characterization of approachable points}\label{sec:finite-laws}

On a probability space we also write $I_\Phi(Z)=\E\Phi(|Z|)$ and
$\|Z\|_F=\inf\{r>0:\E\Phi(|Z|/r)\leq r\}$. Both quantities depend only
on the distribution of $Z$. Throughout this section $\Phi$ denotes the
continuous nondecreasing function chosen in Lemma~\ref{lem:representative}, and $X=L^\Phi(0,1)$ has the modular topology.

\begin{theorem}\label{thm:laws}
The following assertions are equivalent.
\begin{enumerate}
\item[(1)] $\chi_{(0,1)}\in\mathcal A(X)$.
\item[(2)] There exist constants $0<c\leq 1$, $H>0$, and a sequence of
nonnegative random variables $Y_m$, each taking only finitely many values,
such that
\begin{equation}
  \E Y_m=c,\quad \|Y_m\|_F\longrightarrow 0,\quad
  \E\Phi(aY_m)\leq Ha\quad (m\geq 1,\ 0<a\leq 1).
  \label{eq:finite-laws}
\end{equation}
\item[(3)] $\chi_{(0,1)}\in\mathcal N(E_\Phi)$, with the topology induced on $E_\Phi$.
\end{enumerate}
\end{theorem}

To use the laws in Theorem~\ref{thm:laws}, we need one bounded set containing all
convex combinations of their independent copies. The common constant in
\eqref{eq:finite-laws} gives a modular estimate at every sufficiently small scale.

\begin{lemma}\label{lem:rescale}
Let $Y_i$ be independent nonnegative random variables on a probability
space, with $\E Y_i\leq 1$ and $\E\Phi(aY_i)\leq Ha$ for every
$0<a\leq 1$ and every $i$. For every finitely supported family
$\lambda_i\geq 0$ with $\sum_i\lambda_i\leq 1$ and every $0<d\leq 1$,
\begin{equation}
  \E\Phi\left(\frac d8\sum_i\lambda_iY_i\right)\leq\Phi(d)+3Hd.
  \label{eq:independent-sums}
\end{equation}
The closed convex hull of $\{0\}\cup\{Y_i\}$ is topologically bounded in the
Mazur--Orlicz $F$-norm topology.
\end{lemma}

\emph{Proof.} Let $\Psi$ be continuous and nondecreasing, with
$\Psi(0)=0$ and $\Psi(t)>0$ for $t>0$, and assume initially only
\begin{equation}
  \E\bigl[\Psi(aY_i)\chi_{\{aY_i>1\}}\bigr]\leq Ma
  \qquad (0<a\leq 1).
  \label{eq:tail-hypothesis}
\end{equation}
We shall prove
\begin{equation}
  \E\Psi\left(\frac18\sum_i\lambda_iY_i\right)\leq\Psi(1)+3M.
  \label{eq:one-scale-sum}
\end{equation}
Fix the convex coefficients and keep only positive coefficients. Write
\[
  V=\sum_i\lambda_iY_i\chi_{\{\lambda_iY_i>1\}},\qquad
  S=\sum_i\lambda_iY_i\chi_{\{\lambda_iY_i\leq 1\}}.
\]
We first estimate $V$. For each outcome put $I=\{i:\lambda_iY_i>1\}$ and
$N=|I|$. When $N>0$, put $a_i=\min\{N\lambda_i,1\}$ for $i\in I$. The terms
with $N\lambda_i\leq 1$ contribute at most one to
$\sum_{i\in I}\lambda_i/a_i$, and the others contribute at most
$\sum_i\lambda_i\leq 1$. Hence $V\leq 2\max_{i\in I}a_iY_i$, so
\[
  \Psi(V/2)\leq\sum_i\Psi\bigl(\min\{N\lambda_i,1\}Y_i\bigr)
  \chi_{\{\lambda_iY_i>1\}}.
\]
Both sides are zero when $N=0$. Let $N_i=\sum_{j\neq i}\chi_{\{\lambda_jY_j>1\}}$.
On the event $\{\lambda_iY_i>1\}$ we have $N=1+N_i$. Condition on $N_i$,
which is independent of $Y_i$, and set $a=\min\{(1+N_i)\lambda_i,1\}$. Since
$a\geq\lambda_i$, \eqref{eq:tail-hypothesis} bounds the conditional expectation of the $i$th
summand by $Ma\leq M\lambda_i(1+N_i)$. Markov's inequality gives
$\E N_i\leq\sum_j\lambda_j\E Y_j\leq 1$. Taking expectations and summing
yields
\[
  \E\Psi(V/2)\leq 2M.
\]

For $S$, take $t\geq 2$ and put $F(t)=\Psi(t)/t$ for the duration of this
proof. With $x\geq t$, apply \eqref{eq:tail-hypothesis} at $a=t/x$ to obtain
\[
  \Prob\{Y_i>x\}\leq\frac{M}{xF(t)}.
\]
Let $Z_i=\lambda_iY_i\chi_{\{t<Y_i\leq 1/\lambda_i\}}$ and
$T=\sum_i Z_i$. Then $S\leq t+T$, $0\leq Z_i\leq 1$ and $\E T\leq 1$. For
$B>t$,
\begin{align*}
  \E\bigl[Y_i^2\chi_{\{t<Y_i\leq B\}}\bigr]
  &\leq t^2\Prob\{Y_i>t\}+2\int_t^B x\Prob\{Y_i>x\}\,dx\\
  &\leq\frac{M}{F(t)}\bigl(t+2(B-t)\bigr)\leq\frac{2MB}{F(t)}.
\end{align*}
Taking $B=1/\lambda_i$ when $1/\lambda_i>t$, and using $Z_i=0$ otherwise,
yields by independence
\begin{equation}
  \operatorname{Var}(T)\leq\sum_i\E Z_i^2\leq 2M/F(t).
  \label{eq:truncated-variance}
\end{equation}

Fix an order on the finitely many indices and let $\tau$ be the first index
at which a partial sum of the $Z_i$ is at least two. If no such index
exists, put $\tau=\infty$. Then
$\Prob\{\tau<\infty\}=\Prob\{T\geq 2\}\leq 2M/F(t)$ by \eqref{eq:truncated-variance}. On
$\{\tau=j\}$ the partial sum is smaller than three. The remaining sum is
independent of $\{\tau=j\}$. The bound $e^z\leq 1+(e-1)z$ on $[0,1]$ gives
\[
  \E\exp\left(\sum_{i>j}Z_i\right)
  \leq\exp\left((e-1)\sum_{i>j}\E Z_i\right)\leq e^{e-1}.
\]
For $s\geq 3$, the event $T\geq s$ requires $\tau<\infty$. Applying the last
estimate when $\tau=j$ and summing over $j$ gives
\[
  \Prob\{T\geq s\}\leq\frac{2M}{F(t)}e^{e+2-s}.
\]
Partition $\{S\geq 6\}$ into $A_k=\{6k\leq S<6(k+1)\}$ for $k\geq 1$. Use
$t=k+1$ on $A_k$, so $T\geq 5k-1$. For $k=1$, \eqref{eq:truncated-variance} with $\E T\leq 1$ gives
$\Psi(2)\Prob(A_1)\leq 4M/9$. For $k\geq 2$, the exponential bound gives
\[
  \Psi(k+1)\Prob(A_k)\leq 2Me^{e+3}(k+1)e^{-5k}.
\]
Writing $q=e^{-5}$, we have
$\sum_{k\geq 2}(k+1)q^k=q^2(3-2q)/(1-q)^2$, hence
$2e^{e+3}\sum_{k\geq 2}(k+1)e^{-5k}<1/4$. The part with $S<6$ contributes at
most $\Psi(1)$. Therefore $\E\Psi(S/6)\leq\Psi(1)+M$. Since
$(S+V)/8\leq\max\{S/6,V/2\}$, the two estimates prove \eqref{eq:one-scale-sum}.

Finally set $\Psi(u)=\Phi(du)$. Hypothesis \eqref{eq:tail-hypothesis} holds with $M=Hd$, since
$da\leq 1$. Inequality \eqref{eq:one-scale-sum} now gives \eqref{eq:independent-sums}. For $\varepsilon>0$ choose $d$
so that its right-hand side is at most $\varepsilon$. For every
$|b|\leq d\varepsilon/8$ we obtain
$\|b\sum_i\lambda_iY_i\|_F\leq\varepsilon$. This is uniform over all the
convex coefficients, and the bound holds on the closure by continuity. It
proves the asserted topological boundedness. \qed

We first derive the finite distributions in (2) from the absolutely convex
hulls in (1).

\emph{Proof of Theorem~\ref{thm:laws}.} Choose finite rows
$x_{m,1},\ldots,x_{m,n_m}$ at tolerances $t_m\downarrow 0$, with all their
absolutely convex hulls in the same bounded set. We first replace these by
rows of bounded simple functions with mean $\chi_{(0,1)}$. Fix one row and
suppress $m$. For $A=\{\max_i|x_i|\leq T\}$ choose $T$ so large that
$\|\chi_{A^c}\|_F<t_m$, using \eqref{eq:all-scale}. Multiplication by $\chi_A$
gives a mean error smaller than $2t_m$, since
\[
  \chi_{(0,1)}-\chi_A\frac1n\sum_i x_i
  =\chi_{A^c}+\chi_A\left(\chi_{(0,1)}-\frac1n\sum_i x_i\right).
\]
For each signed linear combination the truncated combination is dominated
in absolute value by the original one. Hence boundedness of the union of the
absolutely convex hulls is preserved. Approximate each truncated function
uniformly by a simple function $s_i$ with error at most $e_m$, where
$e_m\downarrow 0$ and $\|e_m\chi_{(0,1)}\|_F<t_m$. For
$d_m=\chi_{(0,1)}-n^{-1}\sum_i s_i$ put $h_i=s_i+d_m$. Then
$\|d_m\|_F<3t_m$, $n^{-1}\sum_i h_i=\chi_{(0,1)}$ and
$\max_i\|h_i\|_F<5t_m$. The absolutely convex combinations of the uniform
approximation errors are bounded in modulus by $e_1$, and the translations
contribute elements of $\{\alpha d_m:|\alpha|\leq 1,\ m\geq 1\}$. Each of
the two sets of errors is topologically bounded. Denote by $V$ the bounded
union of the resulting absolutely convex hulls. Choose $B>1$ such that
$\sup_{v\in V}\|v/B\|_F<1$. Then $I_\Phi(v/B)\leq 1$ for $v\in V$.

Set $g_i=n h_i^+\big/\sum_j h_j^+$. Since $\sum_j h_j=n$, the denominator is
at least $n$. Thus $0\leq g_i\leq|h_i|$ and $\sum_i g_i=n$. Let $J$ be
uniform on $\{1,\ldots,n\}$, independent of the point of $(0,1)$, and set
$Z=g_J$. This has finite support, $\E Z=1$ and
$\E\Phi(Z/(6t_m))\leq 6t_m$. We claim that
\begin{equation}
  \E\bigl[\Phi(aZ/B)\chi_{\{aZ>1\}}\bigr]
  \leq\frac{4a}{1-e^{-1}}\qquad (0<a\leq 1).
  \label{eq:sampled-tail}
\end{equation}
Let $k=\lfloor 1/a\rfloor$ and draw $k$ independent uniform indices
$J_1,\ldots,J_k$, with replacement. Choose independent symmetric signs
$\epsilon_r$, independent also of the indices and the point of the interval.
Since $ka\leq 1$, each $a\sum_r\epsilon_r h_{J_r}$ belongs to $V$ for every
fixed choice of the indices and signs. At a fixed point of the interval,
condition on the indices and all signs except one corresponding to the
largest $|h_{J_r}|$. At least one choice of the remaining sign makes the
absolute value of the signed sum at least $a\max_r|h_{J_r}|$. Consequently
\[
  \E\Phi\left(B^{-1}\max_{r\leq k}a g_{J_r}\right)\leq 2.
\]
At that point let $p=n^{-1}\#\{i:ag_i>1\}$. The identity $\sum_i g_i=n$
gives $p\leq a$, so $kp\leq 1$. For $p>0$, the probability of at least one
index with $ag_{J_r}>1$ is $1-(1-p)^k\geq(1-e^{-1})kp$. Conditional on this
event, the first such index is uniform on the set of these indices. The
maximum dominates the first selected value. Its expected modular at
scale $1/B$ is therefore at least
\[
  (1-e^{-1})\frac{k}{n}\sum_{i:ag_i>1}\Phi(ag_i/B).
\]
This inequality also holds when $p=0$. Integrating it and using
$1/k\leq 2a$ proves \eqref{eq:sampled-tail}. On $\{aZ\leq 1\}$ the integrand $\Phi(aZ/B)$
equals $aZ/B$. Taking $Y_m=Z/B$, $c=1/B$, and
$H=c+4/(1-e^{-1})$ gives \eqref{eq:finite-laws}. Each law can be realized on $(0,1)$ by
taking its values on intervals of the required lengths.

We give the passage from the finite distributions to a needlepoint in
$E_\Phi$.

\begin{lemma}\label{lem:fspace-needle}
Let $Y_{m,j}$, $m,j\geq 1$, be nonnegative bounded random variables on
$(0,1)$. For each $m$, the variables in row $m$ are independent with the
same law as $Y_m$. Suppose that $\E Y_m=c>0$, $\|Y_m\|_F\to 0$, and
\[
  K=\closedco\bigl(\{0\}\cup\{Y_{m,j}:m,j\geq 1\}\bigr)
\]
is topologically bounded. Then the constant $c\chi_{(0,1)}$ is a needlepoint
of $E_\Phi$.
\end{lemma}

\emph{Proof.} We use the choice of coefficient intervals in the proof of
\cite[Theorem 3.5]{KP80}, with the $F$-norm $\|\cdot\|_F$. Fix $\varepsilon>0$ and choose
$\varepsilon_k>0$ with $\sum_{k\geq 1}\varepsilon_k<\varepsilon/12$. Start
with $a_1=1/2$. Given $a_k>0$, choose $m_k$ so that
$\|Y_{m_k}\|_F\leq a_k\varepsilon_k$. Put
$L_k=\|Y_{m_k}-c\|_\infty$ and $v_k^2=\E|Y_{m_k}-c|^2$. For nonnegative,
finitely supported $d_j$ with $\sum_j d_j\leq 1$, the centered sum
$U=\sum_j d_j(Y_{m_k,j}-c)$ satisfies
\[
  |U|\leq L_k,\qquad \E U=0,\qquad
  \E|U|^2=v_k^2\sum_j d_j^2\leq v_k^2\max_j d_j.
\]
Take $z>0$ so small that $\Phi(z)<\varepsilon_k/2$. Chebyshev's inequality
gives
\[
  I_\Phi(U/\varepsilon_k)
  \leq\Phi(z)+\Phi(L_k/\varepsilon_k)
  \frac{v_k^2\max_j d_j}{z^2\varepsilon_k^2}.
\]
Thus we may choose $0<a_{k+1}<a_k/2$ so that $\|U\|_F\leq\varepsilon_k$
whenever $\max_j d_j\leq a_{k+1}$.

Choose $N$ with $\sup_{x\in K}\|x/N\|_F<\varepsilon/4$. Define
\[
  w_j=\frac1N\sum_{k=1}^N Y_{m_k,j}.
\]
Contraction and the triangle inequality give
$\|w_j\|_F\leq\sum_{k=1}^N a_k\varepsilon_k<\varepsilon/12$. Consider any
nonnegative, finitely supported $\lambda_j$ with $\sum_j\lambda_j=1$. At
each $k\leq N$ split its support into
\[
  A_k=\{j:\lambda_j\geq a_k\},\quad
  B_k=\{j:\lambda_j\leq a_{k+1}\},\quad
  C_k=\{j:a_{k+1}<\lambda_j<a_k\}.
\]
There are at most $1/a_k$ indices in $A_k$, so
$\|\sum_{A_k}\lambda_jY_{m_k,j}\|_F\leq\varepsilon_k$. The choice of
$a_{k+1}$ gives
$\|\sum_{B_k}\lambda_j(Y_{m_k,j}-c)\|_F\leq\varepsilon_k$. The coefficient
intervals defining $C_k$ are disjoint as $k$ varies. Thus
$\sum_{k\leq N}\sum_{j\in C_k}\lambda_j\leq 1$ and
$\sum_{k\leq N}\sum_{j\in C_k}\lambda_jY_{m_k,j}\in K$. If
$\theta=N^{-1}\sum_{k\leq N}\sum_{j\in B_k}\lambda_j$, then
$0\leq\theta\leq 1$ and
\[
  \left\|\sum_j\lambda_jw_j-\theta c\chi_{(0,1)}\right\|_F
  \leq 2\sum_{k\leq N}\varepsilon_k+\sup_{v\in K}\|v/N\|_F
  <\varepsilon/2.
\]
For fixed $k\leq N$ the arithmetic means of $Y_{m_k,j}$ converge in
probability to $c$ by the variance estimate, and are uniformly bounded by
a constant depending on the row. They therefore converge in $\|\cdot\|_F$ by
\eqref{eq:all-scale}. The arithmetic means of the $w_j$ converge to $c\chi_{(0,1)}$ as well,
since $N$ is fixed. For some $n$ we have
$\|c\chi_{(0,1)}-n^{-1}\sum_{j\leq n}w_j\|_F<\varepsilon/4$. Add this mean
error to each of $w_1,\ldots,w_n$. The translated family has exact mean
$c\chi_{(0,1)}$, its members have $\|\cdot\|_F$ smaller than $\varepsilon$, and its
convex hull lies within $\varepsilon$ of $[0,c\chi_{(0,1)}]$. All the
functions involved are bounded and belong to $E_\Phi$. This proves the lemma in
the induced topology of $E_\Phi$. \qed

For (2)$\Rightarrow$(3),
realize an array of mutually independent random variables $Y_{m,j}$ on
$(0,1)$ with the row laws given by \eqref{eq:finite-laws}. The existence of countably many
independent random variables on the Lebesgue interval with these prescribed
laws follows from \cite[Theorem 3.19]{Kal02}. The inequality $c\leq 1$ allows us
to apply Lemma~\ref{lem:rescale} to arbitrary finite subfamilies of this entire array. The
resulting closed convex hull with zero is topologically bounded. Lemma~\ref{lem:fspace-needle}
gives $c\chi_{(0,1)}\in\mathcal N(E_\Phi)$. The continuous linear map
$x\mapsto x/c$ on $E_\Phi$ now shows $\chi_{(0,1)}\in\mathcal N(E_\Phi)$. For
(3)$\Rightarrow$(1), the inclusion $E_\Phi\to X$ preserves needlepoints and
$\mathcal N(X)\subset\mathcal A(X)$. \qed

To pass from constants to all of $X$, define
\[
  \theta(x)=\inf\{a>0:I_\Phi(x/a)<\infty\}.
\]
For nondecreasing $\Phi$ and $\alpha,\beta\geq0$ with $\alpha+\beta=1$,
$\Phi(\alpha u+\beta v)\leq\Phi(\max\{u,v\})\leq\Phi(u)+\Phi(v)$; hence
$I_\Phi((x+y)/(a+b))\leq I_\Phi(x/a)+I_\Phi(y/b)$ and $\theta$ is subadditive.
It is absolutely homogeneous and $0\leq\theta\leq\|\cdot\|_F$. Hence it is a continuous
seminorm whose kernel is $E_\Phi$, and $E_\Phi$ is closed. Applying $\theta$ to the
two approximation conditions of an approachable point gives, for every
$\varepsilon>0$,
\[
  \theta(x)\leq\theta\left(x-\frac1n\sum_i x_i\right)
  +\frac1n\sum_i\theta(x_i)<2\varepsilon.
\]
We conclude that $\mathcal A(X)\subset E_\Phi$.

If $x\neq 0$ is approachable, choose $r>0$ with $A=\{|x|\geq r\}$ of
measure $p>0$. Multiplication by $\chi_A/x$, taken to be zero outside
$A$, is continuous and takes $x$ to $\chi_A$. To send $\chi_A$ to
$\chi_{(0,1)}$, take a Borel representative for $A$, and let $\tau$ be the
quantile map of the probability measure $p^{-1}\chi_A\,d\mu$ on
$(0,1)$. Pullback by $\tau$ maps null classes to null classes. For $a>0$,
\[
  I_\Phi(a\,y\circ\tau)
  =p^{-1}\int_A\Phi(a|y(t)|)\,d\mu(t)
  \leq p^{-1}I_\Phi(ay).
\]
It follows from \eqref{eq:all-scale} that this pullback is continuous. We have obtained
$\chi_{(0,1)}\in\mathcal A(X)$ from $x\in\mathcal A(X)\setminus\{0\}$. The reduction to the constant function by a continuous linear operator
is also used in \cite[proof of Theorem 3, p.\ 43]{Ka80}.

By Theorem~\ref{thm:laws} we then have $\chi_{(0,1)}\in\mathcal N(E_\Phi)$. Multiplication by
any bounded function is continuous in $E_\Phi$ and maps $\chi_{(0,1)}$ to that
function. Bounded functions are dense in $E_\Phi$ by truncation and dominated
convergence at each scale. For continuous unbounded Orlicz functions, this
description of $E_\Phi$ as the closure of the bounded functions on a finite
interval already appears in \cite[\S\S4.2--4.4]{Ma60}. Kalton gives the general
form used here in \cite[pp.\ 126--127]{Ka77}. Thus $\mathcal N(E_\Phi)=E_\Phi$ by
closedness of the needlepoint set. In view of
$\mathcal N(E_\Phi)\subset\mathcal N(X)\subset\mathcal A(X)\subset E_\Phi$, we obtain
$\mathcal N(X)=\mathcal A(X)=E_\Phi$. If instead $\chi_{(0,1)}\notin\mathcal A(X)$,
the same level-set argument gives $\mathcal N(X)=\mathcal A(X)=\{0\}$. Since
the inclusion $E_\Phi\to X$ also preserves approachability,
$\mathcal N(E_\Phi)=\mathcal A(E_\Phi)$ equals the same set in each case.

It remains to identify the compact convex sets in the two cases. In the
first case $E_\Phi$ is a nonzero complete separable $F$-space and every point of
$E_\Phi$ is a needlepoint with respect to $E_\Phi$. For separability, truncate and
approximate by bounded simple functions. Their level sets can be
approximated in measure by finite unions of dyadic intervals. Approximating
the finitely many coefficients by rationals gives a countable family dense
at every scale, by \eqref{eq:all-scale}. Theorem 9.4 of \cite{KPR84} now supplies a nonempty
compact convex subset of $E_\Phi$ without extreme points. In fact, the theorem
of Nguyen To Nhu and Le Hoang Tri \cite[Theorem 2-1]{NT94} allows this compact
set to be chosen homeomorphic to the Hilbert cube $[0,1]^{\mathbb N}$. It
has the same topology and the same extreme points in $X$. If it were locally
convex in the relative topology, the embedding theorem of \cite[p.\ 453]{Ro78}
would give an affine homeomorphism onto a compact convex subset of a locally
convex space, contrary to the Krein--Milman theorem \cite{KM40}. In the other
case, zero is the only approachable point of $X$, so \cite[Theorem 1]{Ka80}
embeds every compact convex subset of $X$ affinely and homeomorphically in
a locally convex space. This gives both relative local convexity and the
compact Krein--Milman property. This proves Theorem~\ref{thm:modular}, and the equality
of the two topologies in Lemma~\ref{lem:gauge} gives Theorem~\ref{thm:lux}.

\emph{Orlicz functions bounded above by a linear function at
infinity.} Under the standing assumption $\varphi\in\Delta_2$, Kalton and
Peck proved
that the unit ball of $L^1(0,1)$ contains a compact convex set without
extreme points in the topology of $L^\varphi(0,1)$ when
$\varphi(u)/u$ has lower limit zero and finite upper limit at infinity
\cite[Section 4, Example 4.3]{KP80}. For continuous increasing concave
functions $\varphi$, positive away from zero and satisfying $\varphi(u)/u\to 0$, De
Pascale, Trombetta and Weber obtained strong needlepoints in
$L^\varphi(0,1)$ \cite[Definition 1.1 and proof of Corollary 1.8]{DTW93}. We
drop the assumption $\varphi\in\Delta_2$ and identify the subspace that
consists of needlepoints. It
can be a proper subspace of $L^\varphi(0,1)$. We also characterize exactly
when the distributions in Theorem~\ref{thm:laws} may have just one nonzero value.

\begin{proposition}\label{prop:linear-growth}
Suppose that $\limsup_{u\to\infty}\varphi(u)/u<\infty$. Then the inclusion
of $L^1(0,1)$ in $E_\varphi$ is continuous, and the closure of $L^1(0,1)$
in $X=L^\varphi(0,1)$ is $E_\varphi$. The following alternatives hold.
\begin{enumerate}
\item[(1)] If $\liminf_{u\to\infty}\varphi(u)/u=0$, then
$\mathcal A(X)=\mathcal N(X)=E_\varphi$, and $E_\varphi$ contains a
nonempty compact convex set without extreme points.
\item[(2)] If $\liminf_{u\to\infty}\varphi(u)/u>0$, then
$X=E_\varphi=L^1(0,1)$ with the usual norm topology, and
$\mathcal A(X)=\mathcal N(X)=\{0\}$.
\end{enumerate}
More precisely, for a continuous representative $\Phi$ as in Section~\ref{sec:prelim}, the
laws in Theorem~\ref{thm:laws} can all be supported on zero and one positive value if
and only if
\[
  \limsup_{u\to\infty}\Phi(u)/u<\infty
  \quad\text{and}\quad
  \liminf_{u\to\infty}\Phi(u)/u=0.
\]
This last assertion holds without the initial upper bound on $\varphi$.
\end{proposition}

\emph{Proof.} For the representative constructed in Section~\ref{sec:prelim}, the
inequalities
\[
  \frac{\varphi(u)}u\leq\frac{G(u)}u\leq 2\frac{\varphi(2u)}{2u}
\]
show that finiteness of the upper limit and vanishing of the lower limit
are each preserved in both directions. Rescaling the argument and dividing
by $G(b)>0$ preserve these same properties, and the change on $[0,1]$ does
not affect them. Put $F(u)=\Phi(u)/u$. Since $F=1$ on $(0,1]$ and is
continuous on $[1,\infty)$, a finite upper limit gives
$M=\sup_{u>0}F(u)<\infty$. For $f\in L^1(0,1)$ and $a>0$ we have
$I_\Phi(af)\leq Ma\|f\|_1$. This proves $L^1(0,1)\subset E_\varphi$ and
continuity of the inclusion by \eqref{eq:all-scale}. The closure is $E_\varphi$, since
$E_\varphi$ is closed in $X$ and contains the bounded functions as a dense
subspace.

Suppose that $\liminf_{u\to\infty}F(u)=0$. Choose $s_m\geq m^2$ with
$F(s_m)\leq m^{-3}$ and set $u_m=s_m/m$. Take a set $A_m$ of measure
$1/u_m$ and put $Y_m=u_m\chi_{A_m}$. Then $\E Y_m=1$ and, for
$0<a\leq 1$,
\[
  \E\Phi(aY_m)=aF(au_m)\leq Ma.
\]
For each fixed $b>0$ we have $bu_m\leq s_m$ when $m\geq b$, hence
\[
  I_\Phi(bY_m)=\frac{\Phi(bu_m)}{u_m}\leq mF(s_m)\leq m^{-2}.
\]
Thus $\|Y_m\|_F\to 0$ by \eqref{eq:all-scale}. Theorem~\ref{thm:laws} and the preceding
classification give (1).

Conversely, suppose the laws of Theorem~\ref{thm:laws} have just one positive value
$u_m$ besides zero. Their common mean $c>0$ forces the mass at $u_m$ to be
$c/u_m$, so $u_m\geq c$. Since $Y_m\to 0$ in measure, $c/u_m\to 0$ and
consequently $u_m\to\infty$. Their modular at scale one is
$cF(u_m)\to 0$. For any $v>0$, choose an $m$ with $u_m\geq v$ and take
$a=v/u_m$ in \eqref{eq:finite-laws}. The resulting inequality is $caF(v)\leq Ha$, so
$F(v)\leq H/c$ on the entire positive half-line. Together with the
construction above, this proves the assertion about the supports.

Finally, if the lower limit is positive, then $\inf_{u>0}F(u)>0$ by
positivity and continuity and the identity $F=1$ near zero. Thus
$m_0u\leq\Phi(u)\leq Mu$ for some $m_0>0$ and all $u\geq 0$. Integrating at
arbitrary scales gives $X=E_\varphi=L^1(0,1)$ with its norm topology. In a
normed space the norm of an arithmetic mean is at most the maximum norm of
its summands, so only zero is approachable. This proves (2). \qed

For example, fix $0<p<1$ and put $r_n=2^{n(n+1)}$ for $n\geq 0$. Set
$\Phi(t)=t$ on $[0,1]$, define $\Phi$ on $[1,\infty)$ by
\[
  \Phi(t)=\max\{r_n^{1-p}t^p,\,2t-r_{n+1}\}
  \qquad (r_n\leq t\leq r_{n+1}),
\]
and extend $\Phi$ evenly to $\mathbb R$.
The two formulas at each $r_n$ give $\Phi(r_n)=r_n$, and $\Phi(t)\leq t$.
The ratio $\Phi(t)/t^p$ is nondecreasing: on each interval it is the
maximum of a constant and $2t^{1-p}-r_{n+1}t^{-p}$, an increasing function.
Thus $\Phi(2^{1/p}t)\geq 2\Phi(t)$ and $\Phi$ is $b$-Orlicz. For
$t_n=r_{n+1}/2$ and $q_n=t_n/r_n=2^{2n+1}$ we have
$\Phi(t_n)/t_n=q_n^{p-1}\to 0$. Proposition~\ref{prop:linear-growth} therefore applies with
two-point laws. To see that $E_\Phi$ is proper, take disjoint sets $C_n$ of
measure $1/r_{n+1}$, possible since $r_{n+1}\geq 4^{n+1}$. For
$f=\sum_{n\geq 0}t_n\chi_{C_n}$ one has
\[
  I_\Phi(f)=\frac12\sum_{n\geq 0}q_n^{p-1}<\infty,\qquad
  I_\Phi(2f)=\sum_{n\geq 0}\Phi(r_{n+1})/r_{n+1}=\infty.
\]

In the example below, $\Phi(u)/u$ is unbounded above at infinity. Proposition~\ref{prop:linear-growth} therefore excludes two-point laws with a single positive value. We
construct the required finite distributions explicitly.

We use the Matuszewska--Orlicz indices at infinity $\alpha_\Phi$ and
$\beta_\Phi$, defined in Section~\ref{sec:prelim} (Definition~\ref{def:indices}) in the sense of
\cite[Definition~1.7]{KZ19}.

\begin{proposition}\label{prop:proper-example}
For each $0<p<1$ there is a $b$-Orlicz function $\Phi$, with
$b=2^{1/p}$, for which
\[
  \mathcal A(L^\Phi)=\mathcal N(L^\Phi)=\mathcal A(E_\Phi)
  =\mathcal N(E_\Phi)=E_\Phi\subsetneq L^\Phi,\qquad
  L^\Phi=L^\Phi(0,1).
\]
The subspace $E_\Phi$ contains a nonempty compact convex set without
extreme points. The function $\Phi$ may be chosen so that
\[
  \liminf_{u\to\infty}\frac{\Phi(u)}u=0,\qquad
  \limsup_{u\to\infty}\frac{\Phi(u)}u=\infty,
\]
\[
  \liminf_{B\to\infty}\frac1{\ln B}\int_1^B\frac{\Phi(u)}{u^2}\,du
  =\frac1{12}.
\]
Its lower and upper dilation indices at infinity are $p$ and $\infty$,
respectively. A Poisson random variable of mean one does not belong to
$L^\Phi$.
\end{proposition}

\emph{Proof.} Set $c=1-p$ and $u_1=1/c$. Define $\sigma(s)=0$ on
$(-\infty,0]$ and $\sigma(s)=-cs$ on $[0,u_1]$. Having defined $\sigma$ up
to $u_n$, put
\[
  J_n=\int_0^{u_n}(e^{\sigma(t)}-1)_+\,dt,\qquad
  L_n=12(1+J_n)+4n+4,
\]
\[
  v_n=u_n+L_n,\qquad w_n=v_n+1,\qquad h_n=e^{2w_n},
\]
\[
  u_{n+1}=w_n+(h_n+n+1)/c.
\]
Extend $\sigma$ continuously by setting
\[
  \sigma(s)=
  \begin{cases}
    -n, & u_n\leq s\leq v_n,\\
    -n+(h_n+n)(s-v_n), & v_n\leq s\leq w_n,\\
    h_n-c(s-w_n), & w_n\leq s\leq u_{n+1}.
  \end{cases}
\]
The final value is $-(n+1)$, as required at the next step. Let
$\Phi(0)=0$ and $\Phi(u)=|u|e^{\sigma(\ln|u|)}$ for $u\neq 0$. This function
equals $|u|$ on $[-1,1]$. All slopes of $s+\sigma(s)$ are at least $p$.
Hence $\Phi(au)\geq a^p\Phi(u)$ for $a\geq 1$, $u>0$, which gives the
asserted $b$. On $[w_n,u_{n+1}]$ the slope equals $p$ and the lengths tend
to infinity. A lower dilation bound $\Phi(au)\geq Ca^r\Phi(u)$ for all
$a\geq 1$ and all sufficiently large $u$ is therefore impossible for
$r>p$. Further,
\[
  \frac{\Phi(2e^{v_n})}{\Phi(e^{v_n})}=2^{1+h_n+n}\longrightarrow\infty,
\]
so there is no finite upper dilation index. Since $\sigma(u_n)=-n$ and
$\sigma(w_n)=h_n$, the asserted extreme limits of $\Phi(u)/u$ follow.

We verify the hypotheses of Theorem~\ref{thm:laws} using finite distributions. Put
$k_n=\lfloor L_n/2\rfloor$ and $z_j=u_n+L_n/4+j$ for $0\leq j<k_n$. Give
$Y_n$ the value $e^{z_j}$ with probability $e^{-z_j}/k_n$ and give it the
value zero with the remaining probability. Thus $\E Y_n=1$. For
$0<a\leq 1$,
\[
  \frac{\E\Phi(aY_n)}a
  =\frac1{k_n}\sum_{j=0}^{k_n-1}e^{\sigma(z_j+\ln a)}.
\]
Since $\sigma'\geq-c$ almost everywhere, $e^{\sigma(r)}\leq
e^c e^{\sigma(t)}$ for $r\leq t\leq r+1$. The unit intervals starting at
$z_j+\ln a$ have disjoint interiors and lie below $v_n$. The integral of
$(e^\sigma-1)_+$ on $(-\infty,v_n]$ equals $J_n$. Consequently the last
mean is at most $e^c(1+J_n/k_n)\leq 2e$, uniformly in $a$ and $n$. For
fixed $a>0$, all $z_j+\ln a$ belong to $[u_n,v_n]$ once
$L_n/4>|\ln a|$. It follows that $I_\Phi(aY_n)=ae^{-n}$ for all
sufficiently large $n$. By \eqref{eq:all-scale}, $\|Y_n\|_F\to 0$. Theorems~\ref{thm:laws} and~\ref{thm:lux} give the stated equalities and the compact convex set in $E_\Phi$.

To see that $E_\Phi$ is proper, put $d_n=e^{v_n}$ and choose disjoint measurable sets
$C_n\subset(0,1)$ with measure $2^{-n}/\Phi(d_n)$, possible since
$\Phi(d_n)>1$. Then $f=\sum_{n\geq 1}d_n\chi_{C_n}$ satisfies
\[
  I_\Phi(f)=\sum_{n\geq 1}2^{-n}=1,\qquad
  I_\Phi(2f)=\sum_{n\geq 1}2^{-n}2^{1+h_n+n}=\infty.
\]
For any $a>0$ let $m_n=\lceil ae^{w_n}\rceil$. Monotonicity gives
\[
  \ln\frac{\Phi(m_n/a)}{m_n!}
  \geq w_n+e^{2w_n}-m_n\ln m_n\longrightarrow\infty.
\]
Thus $\sum_{m\geq 1}\Phi(m/a)/m!=\infty$ at every scale, as claimed for
the Poisson distribution.

Finally write $H(s)=\int_0^s e^{\sigma(t)}\,dt$. On the descending
interval $[w_n,u_{n+1}]$, integration over $\sigma\geq 0$ gives
$J_{n+1}\geq(e^{h_n}-1-h_n)/c$. Since $u_{n+1}=O(h_n)$ and
$n=O(h_n)$, we have $u_n/J_n\to 0$, $n/J_n\to 0$ and
$v_n/J_n,w_n/J_n\to 12$. The inequalities $J_n\leq H(u_n)\leq J_n+u_n$ and
the identity $H(v_n)=H(u_n)+e^{-n}L_n$ give $H(v_n)/v_n\to 1/12$. For
$u_n\leq s\leq w_n$, $H(s)/s\geq J_n/w_n$. For $w_n\leq s\leq u_{n+1}$,
integration over the preceding interval of length one gives
\[
  \frac{H(s)}s\geq\frac{e^{h_n}-e^{-n}}{(h_n+n)u_{n+1}}\longrightarrow\infty.
\]
This proves $\liminf_{s\to\infty}H(s)/s=1/12$. Substituting $u=e^t$ gives
the stated logarithmic average. \qed

\section{An approachable point which is not a needlepoint}\label{sec:counter}

We turn to the space in Theorem~\ref{thm:counter}. Throughout the construction,
$p=1/2$. In the $k$th summand we shall exhibit a unit vector which is the
arithmetic mean of vectors of norm at most $2^{-k}$.

Fix an integer $S\geq 2$ and put
\[
  K=4096(S+1)^3,\qquad m=16(K+1)S.
\]
On $E=\mathbb R^m$, let $\|\cdot\|_1$ and $\|\cdot\|_p$ denote the usual
$\ell_1$-norm and $p$-norm, respectively. Set
\[
  q(x)=\inf_{b+c=x}\left(\|b\|_1^p+S^{-p}\|c\|_p^p\right)^{1/p},
  \qquad e=m^{-1}(1,\ldots,1).
\]
The infimum is attained since its sublevel sets are compact. The
$p$-triangle inequality follows by adding decompositions. For every
$b\in E$, Cauchy--Schwarz gives
$\|b\|_p=(\sum_j|b_j|^{1/2})^2\leq m\|b\|_1$. Since $S\leq m$, the
$p$-triangle inequality for $m^{-1}\|\cdot\|_p$ gives
\[
  m^{-1}\|x\|_p\leq q(x)\leq\min\{\|x\|_1,S^{-1}\|x\|_p\}.
\]
Thus $q$ is a $p$-norm, continuous in the Euclidean topology, and
$q(e)=1$. We write $\operatorname{dist}_q(x,\mathbb R e)
=\inf_{a\in\mathbb R}q(x-ae)$. This infimum is attained, since
$q(x-ae)^p\geq|a|^p-q(x)^p$.

Suppose $q(x)\leq h$ and take a minimizing decomposition $x=b+c$. Then
$\|b\|_1\leq h$ and $\|c\|_p\leq Sh$. The $j$th largest coordinate
magnitude of $c$ is at most $\|c\|_p j^{-2}$. Retaining the indices $I$ of
its $S$ largest coordinates and using $\sum_{j>S}j^{-2}\leq 1/S$, we
obtain
\begin{equation}
  \|x\|_1\leq(S+1)h,\qquad |I|\leq S,\qquad
  \sum_{j\notin I}|x_j|\leq 2h.
  \label{eq:counter-sparse}
\end{equation}

\begin{lemma}\label{lem:counter-block}
Let $v_1,\ldots,v_n\in E$ and $w_i>0$ satisfy $\sum_i w_i=1$ and
$\sum_i w_i v_i=\beta e$. Put
\[
  t=\sum_i w_i q(v_i),\qquad
  s=\sup_{\sum_i|a_i|\leq 1}
  \operatorname{dist}_q\left(\sum_i a_i v_i,\mathbb R e\right).
\]
Then $|\beta|\leq 64\max\{t,s\}$.
\end{lemma}

\emph{Proof.} We first consider a finitely supported random vector $W$ of
mean $e$. Put $h=1/32$. Suppose that every value of $W$ has $q$-norm at
most $h$ and every convex combination of its values is within $q$-distance
$h$ of $\mathbb R e$. Let $Y$ be the mean of $K$ independent copies of $W$,
and let $\ell(x)=\sum_j x_j$. By \eqref{eq:counter-sparse}, both $\|W\|_2$ and $|\ell(W)|$ are
at most $S+1$. Independence gives
\[
  \E\|Y-e\|_2^2\leq(S+1)^2/K,\qquad
  \E|\ell(Y)-1|^2\leq(S+1)^2/K.
\]
Choose a realization for which
$S\|Y-e\|_2^2+|\ell(Y)-1|^2\leq(S+1)^3/K=1/4096$. For every set $V$ of at
most $S$ indices, Cauchy--Schwarz yields
\[
  \ell(Y)\geq 1-1/64,\qquad
  \sum_{j\in V}|Y_j|\leq S/m+\sqrt S\,\|Y-e\|_2\leq 5/64.
\]
Apply \eqref{eq:counter-sparse} to each sampled value of $W$ and let $U$ be the union of the
resulting sets of indices. Then $|U|\leq KS$ and
$\sum_{j\notin U}|Y_j|\leq 2h$. Our assumption on convex combinations
gives $a\in\mathbb R$ with $q(Y-ae)\leq h$. Another application of \eqref{eq:counter-sparse}
gives $V$, with $|V|\leq S$, such that
$\sum_{j\notin V}|Y_j-a/m|\leq 2h$. As $(K+1)S/m=1/16$, summing over the
complement of $U\cup V$ gives $15|a|/16\leq 4h$. On the other hand,
\[
  1-6/64\leq\ell(Y)-\sum_{j\in V}|Y_j|
  \leq\sum_{j\notin V}|Y_j|\leq|a|+2h.
\]
These two bounds contradict $h=1/32$.

Now suppose $|\beta|>64\max\{t,s\}$. Then $t>0$, for otherwise $v_i=0$
for all $i$. The vectors and weights
\[
  b_i=\frac{v_i}{1+q(v_i)/t},\qquad
  w_i'=\frac12 w_i\bigl(1+q(v_i)/t\bigr)
\]
satisfy $\sum_i w_i'=1$, $\sum_i w_i'b_i=\beta e/2$ and $q(b_i)\leq t$.
Every signed combination of the $b_i$ with $\sum_i|a_i|\leq 1$ is such a
combination of the $v_i$. Hence the random vector taking values
$2b_i/\beta$ with probabilities $w_i'$ has mean $e$, its values have
$q$-norm less than $1/32$, and their convex combinations have $q$-distance
less than $1/32$ from $\mathbb R e$. This contradicts the first part of the
proof. \qed

\begin{lemma}\label{lem:counter-operators}
Let $(B_k)$ be $p$-normed spaces, with $p=1/2$, and let
$T:\ell_1^n\to(\bigoplus_k B_k)_p$ be a bounded linear operator with
coordinate operators $T_k$. Then $\sum_k\|T_k\|\leq 2\|T\|$.
\end{lemma}

\emph{Proof.} Fix finitely many indices $k$ with $s_k=\|T_k\|>0$. By
compactness, choose $v_k$ in the unit ball of $\ell_1^n$ with
$\|T_kv_k\|=s_k$. For positive $a_k$ summing to one and independent
random signs $\varepsilon_k$, the vector $v=\sum_k a_k\varepsilon_kv_k$
belongs to the unit ball. The $p$-triangle inequality gives
$2^p\|y\|^p\leq\|x+y\|^p+\|x-y\|^p$. Holding the other signs fixed and
averaging over $\varepsilon_k$ therefore yields
$\E\|T_kv\|^p\geq 2^{p-1}a_k^p s_k^p$. Summing in $k$ gives
\[
  \|T\|^p\geq 2^{p-1}\sum_k a_k^p s_k^p.
\]
Set $a_k=s_k/\sum_j s_j$. Since $2p=1$, the sum on the right before
multiplication by $2^{p-1}$ equals $(\sum_k s_k)^p$. Taking the supremum
over finite sets proves the assertion. \qed

For $k\geq 1$, take $S_k=2^k$ and denote the corresponding objects constructed above by $E_k,q_k,m_k,e_k$. We also write $e_k$ for the element of
$D=(\bigoplus_k E_k)_p$ supported on the $k$th coordinate with value
$e_k$. Let $N$ be the closed span of the differences $e_k-e_1$ for
$k\geq 2$. Then
\begin{equation}
  N=\left\{(\alpha_k e_k):(\alpha_k)\in\ell_p,\ \sum_k\alpha_k=0\right\}.
  \label{eq:counter-kernel}
\end{equation}
Indeed, the subspace of all $(\alpha_k e_k)_k$ is closed and isometric to
$\ell_p$, and summation on it is continuous because
$\sum_k|\alpha_k|\leq\|(\alpha_k)\|_p$. Conversely, if $\sum_k\alpha_k=0$,
subtract the sum of each finite truncation from its first coordinate. The
resulting zero-sum truncations converge in $\ell_p$.

Put $X=D/N$ with the quotient $p$-norm and $u=e_1+N$. The space $D$ is
complete by coordinatewise completeness and summability of the $p$th
powers. Its quotient by the closed subspace $N$ is also complete. Indeed,
from any Cauchy sequence in the quotient take a subsequence whose
successive differences have summable $p$th powers of their norms. Lift
those differences with at most twice their quotient norms. The series of
the lifts converges in $D$, so the subsequence, and hence the original
sequence, converges in $X$. Separability follows from that of $D$.

All the $e_k$ have the same image $u$. Every lift of $u$ has the form
$(\alpha_k e_k)_k$ with $\sum_k\alpha_k=1$, so its norm is at least one.
Since $q_1(e_1)=1$, we have $\|u\|=1$. The $m_k$ coordinate unit vectors
of $E_k$ have images in $X$ of norm at most $1/S_k$ and arithmetic mean
$u$. Every signed combination of their images has norm at most the
$\ell_1$-norm of its coefficients, because $q_k\leq\|\cdot\|_1$. Thus $u$
and all its scalar multiples are approachable.

Consider the contraction
\[
  Q:X\longrightarrow\left(\bigoplus_k E_k/\mathbb R e_k\right)_p,
  \qquad Q(y+N)=(y_k+\mathbb R e_k)_k.
\]
The formula descends to the quotient, and taking the infimum over lifts
proves the norm bound. By \eqref{eq:counter-kernel}, $\ker Q=\mathbb R u$. Indeed, if
$y_k=\alpha_k e_k$ for all $k$, then $y+N=(\sum_k\alpha_k)u$. For
$x\notin\mathbb R u$, a linear functional on a nonzero finite-dimensional
coordinate of $Qx$ yields a continuous linear functional on $X$ nonzero at
$x$. Every continuous linear functional vanishes at an approachable point
or a needlepoint, since it tends uniformly to zero on the vectors whose
arithmetic means approximate that point. Hence the approachable points of
$X$ are precisely $\mathbb R u$, and every needlepoint belongs to
$\mathbb R u$.

For both remaining assertions we shall use the following consequence of
Lemma~\ref{lem:counter-block}. Let $z_i\in X$ have norms at most $d$, and suppose
$\sum_i w_i z_i=\gamma u$, where $w_i>0$ and $\sum_i w_i=1$. Choose lifts
$y_i=(y_{ik})_k\in D$ with $\|y_i\|_D\leq 2d$. By \eqref{eq:counter-kernel},
\[
  \sum_i w_i y_{ik}=\beta_k e_k,\qquad(\beta_k)\in\ell_p,\qquad
  \sum_k\beta_k=\gamma.
\]
Let $T$ be the operator on $\ell_1^n$ with columns $Qz_i$. Put
$s_k=\|T_k\|$ and $t_k=\sum_i w_i q_k(y_{ik})$. The contractive
inclusion of $\ell_p$ in $\ell_1$ gives
\[
  \sum_k t_k=\sum_i w_i\sum_k q_k(y_{ik})
  \leq\sum_i w_i\|y_i\|_D\leq 2d.
\]
The norm $s_k$ is the supremum of
$\operatorname{dist}_{q_k}(\sum_i a_i y_{ik},\mathbb R e_k)$ for
$\sum_i|a_i|\leq 1$. Lemma~\ref{lem:counter-block} therefore gives, for every set $J$ of
indices,
\begin{equation}
  \sum_{k\in J}|\beta_k|\leq 64\sum_{k\in J}(t_k+s_k),
  \qquad\sum_k t_k\leq 2d.
  \label{eq:counter-lift}
\end{equation}

\emph{Exclusion of nonzero needlepoints.} By homogeneity it suffices to
exclude $u$. Suppose it admits a needlepoint family $x_1,\ldots,x_n$ at
tolerance $h>0$. Set $c=u-n^{-1}\sum_i x_i$ and $z_i=x_i+c$. The $z_i$
have arithmetic mean $u$ and norms at most $2^{1/p}h=4h$. Every convex
combination of the $z_i$ is within $4h$ of $\mathbb R u$, by the segment
condition for the $x_i$. Consequently every convex combination of the
columns $Qz_i$ has norm at most $4h$. For a signed combination of total
variation at most one, let $\alpha,\delta$ be its positive and negative
masses. Its norm is at most
$4h(\alpha^p+\delta^p)^{1/p}\leq 8h$. Thus $\|T\|\leq 8h$ and Lemma~\ref{lem:counter-operators}
gives $\sum_k s_k\leq 16h$. Apply \eqref{eq:counter-lift} with $d=4h$, $w_i=1/n$ and all
indices. It follows that
\[
  1=\left|\sum_k\beta_k\right|\leq 64(8h+16h)=1536h.
\]
Taking $h<1/1536$ proves the assertion. \qed

\emph{Relative local convexity.} Let $C$ be a nonempty compact absolutely
convex subset of $X$, let $B_X$ denote the closed unit ball of $X$, and
put $H_h=\closedco(C\cap hB_X)$. We first show that
$\bigcap_{h>0}H_h=\{0\}$. Every continuous linear functional vanishes on
the intersection, since its supremum in absolute value on
$C\cap hB_X$ tends to zero. The separation by finite-dimensional
coordinates established above thus places every $x$ in the intersection
on $\mathbb R u$. Write $x=\gamma u$.

For each $h>0$, choose $x_i\in C\cap hB_X$ and positive $w_i$ summing to
one such that $c=x-\sum_i w_i x_i$ satisfies $\|c\|\leq h$. We have
$c\in 2C$, so $z_i=x_i+c$ belong to $3C$, have norms at most $4h$, and
satisfy $\sum_i w_i z_i=\gamma u$. With the preceding notation, $T$ maps
the unit ball of $\ell_1^n$ into the fixed compact set $Q(3C)$, by
absolute convexity.

Let $P_{>j}$ be the projection onto the coordinates of index greater than
$j$. Then $\eta_j=\sup_{v\in Q(3C)}\|P_{>j}v\|$ tends to zero. Indeed,
the projections are contractions converging pointwise to zero, so a
finite net for the compact set and the $p$-triangle inequality give
uniform convergence there. Applying Lemma~\ref{lem:counter-operators} to $P_{>j}T$ gives
$\sum_{k>j}s_k\leq 2\eta_j$. By \eqref{eq:counter-lift} with $d=4h$,
\[
  \sum_{k>j}|\beta_k|\leq 64(8h+2\eta_j).
\]
In each of the first $j$ coordinates, $q_k(y_{ik})\leq 8h$, so \eqref{eq:counter-sparse} gives
$|\beta_k|=\|\beta_k e_k\|_1\leq\sum_i w_i\|y_{ik}\|_1\leq 8h(S_k+1)$.
Hence
\[
  |\gamma|\leq 8h\sum_{k\leq j}(S_k+1)+64(8h+2\eta_j).
\]
First let $h\downarrow 0$ at fixed $j$, then let $j\to\infty$. This proves
$\gamma=0$.

The sets $H_h$ are nested compact convex subsets of $C$, each containing
the relative neighborhood $C\cap\{v:\|v\|<h\}$ of zero. By compactness
and the intersection just proved, every relative neighborhood of zero
contains some $H_h$. For a general compact convex set $A$ and $z\in A$,
apply this to $C=\co((A-z)\cup(z-A))$. This is an absolutely convex
compact set, since it is the image of $[0,1]\times(A-z)\times(z-A)$ under
$(t,a,b)\mapsto ta+(1-t)b$. Intersect the convex neighborhoods in $C$
with $A-z$ and translate by $z$. This gives a convex relative
neighborhood base at every point of $A$. The embedding results of Lawson
\cite[Corollary 4.2]{La76} and Roberts \cite[p.\ 453]{Ro78} apply to a base of
convex neighborhoods, including the closed neighborhoods constructed
above. They give an affine homeomorphism of $A$ onto a compact convex
subset of a Hausdorff locally convex space. The classical Krein--Milman
theorem now applies to every nonempty $A$. This completes the proof of
Theorem~\ref{thm:counter}. \qed

Ribe's construction gives a nonlocally convex space with a one-dimensional
uncomplemented subspace and quotient isomorphic to $\ell_1$ \cite[Theorem 1]{Ri79}. Every vector in the kernel of such a nonsplit one-dimensional
extension of a Banach space is a needlepoint \cite[p.\ 1111]{Ka03}. In the
present space the line $\mathbb R u$ contains no nonzero needlepoint,
although the continuous dual annihilates it. Write $J_X$ for the map from
$X$ into its Banach envelope $\widehat X$. When the dual separates
points, the envelope norm on $X$ is the gauge of the convex hull of the
unit ball \cite[p.\ 65]{Ka86}. For general $X$, the envelope seminorm
$x\mapsto\|J_Xx\|$ is the largest seminorm dominated by $\|\cdot\|_X$,
and $J_X(X)$ is dense in $\widehat X$ \cite[\S1]{AAW21}.

\begin{proposition}\label{prop:counter-envelope}
For each $k$, equip $F_k=\mathbb R^{m_k}/\mathbb R e_k$ with the norm
\[
  d_k(v+\mathbb R e_k)
  =S_k^{-1}\inf_{a\in\mathbb R}\|v-ae_k\|_1.
\]
Then $\widehat X$ is isometric to $(\bigoplus_k F_k)_1$, with envelope
map $J_X(y+N)=(y_k+\mathbb R e_k)_k$. In particular,
$\ker J_X=\mathbb R u=\bigcap_{f\in X^*}\ker f$.
\end{proposition}

\emph{Proof.} For $v=b+c$ we have
\[
  S_k^{-1}\|v\|_1
  \leq\|b\|_1+S_k^{-1}\|c\|_p
  \leq\left(\|b\|_1^p+S_k^{-p}\|c\|_p^p\right)^{1/p}.
\]
Thus $S_k^{-1}\|v\|_1\leq q_k(v)$. The map in the statement descends to
$X$ by \eqref{eq:counter-kernel} and is contractive, since
$\sum_k q_k(y_k)\leq\|y\|_D$. Denote this map by $J$ and the envelope
seminorm by $\nu$. We have $\|Jx\|\leq\nu(x)$.

Let $b_{kj}$ be the image in $X$ of the $j$th coordinate unit vector in
$E_k$. Since $u=m_k^{-1}\sum_j b_{kj}$ and $\|b_{kj}\|_X\leq S_k^{-1}$,
we have $\nu(u)\leq S_k^{-1}$ for every $k$, hence $\nu(u)=0$. For
$v\in E_k$, write $[v]$ for its image in $X$. Subtracting $ae_k$ changes
$[v]$ by $au$. Expansion in the coordinate unit vectors therefore gives
\[
  \nu([v])=\nu([v-ae_k])\leq S_k^{-1}\|v-ae_k\|_1
  \qquad (a\in\mathbb R).
\]
Taking the infimum in $a$ and summing over finitely many coordinates
gives $\nu(x)\leq\|Jx\|$ for $x$ with a finitely supported lift. Such
elements are dense, and both seminorms are continuous. Thus
$\nu(x)=\|Jx\|$ throughout $X$. The image contains every finitely
supported element of $(\bigoplus_k F_k)_1$, which proves the isometry
after completion. The kernel is $\mathbb R u$ by \eqref{eq:counter-kernel}, and the last
equality follows from Hahn--Banach. \qed

On the compact segment $[0,u]$ the affine function $tu\mapsto t$ is
continuous. Every continuous affine function on $X$ is constant on this
segment. Thus even on a line segment the continuous affine functions
supplied by relative local convexity need not extend to $X$.

\medskip
\noindent Funding: This research was supported by Heilongjiang Provincial Natural Science Foundation of China (Grant No.~PL2025A002).

\vspace{1.5ex}

\begin{center}
\small
\begin{minipage}{0.92\linewidth}
\scshape
Department of Mathematics, Harbin University of Science and Technology,
Harbin 150080, China, and Department of Mathematics, Harbin University,
Harbin 150086, China\\
\normalfont\itshape
\small Email address: jidong@hrbu.edu.cn
\end{minipage}

\vspace{1.2ex}

\begin{minipage}{0.92\linewidth}
\scshape
Department of Mathematics, Harbin University of Science and Technology,
Harbin 150080, China\\
\normalfont\itshape
\small Email address: cuiya@hrbust.edu.cn
\end{minipage}
\end{center}

\end{document}